\documentclass[12pt]{article}
\usepackage{amsmath}
\usepackage{amsfonts}
\usepackage{amsthm}
\normalbaselines
\theoremstyle{plain}
\newtheorem{theorem}{Theorem}[section]
\newtheorem{proposition}{Proposition}[section]
\newtheorem{corollary}{Corollary}[section]

\theoremstyle{definition}
\newtheorem{definition}{Definition}[section]
\theoremstyle{remark}

\catcode `\@=11
\@addtoreset{equation}{section}

\def\section{\@startsection {section}{1}{\z@}{-3.5ex plus -1ex minus
     -.2ex}{2.3ex plus .2ex}{\normalsize\bf}}
\def\subsection{\@startsection{subsection}{2}{\z@}{-3.25ex plus -1ex minus
 -.2ex}{1.5ex plus .2ex}{\normalsize\bf}}

\def\thebibliography#1{\section*{References\markboth
  {REFERENCES}{REFERENCES}}\list
  {[\arabic{enumi}]}{\settowidth\labelwidth{[#1]}\leftmargin\labelwidth
  \advance\leftmargin\labelsep
  \usecounter{enumi}}
  \def\newblock{\hskip .11em plus .33em minus -.07em}
  \sloppy
  \sfcode`\.=1000\relax}
 
\catcode `\@=12
\newcommand{\C}{\ensuremath{\mathbb{C}}}
\newcommand{\N}{\ensuremath{\mathbb{N}}}

\renewcommand{\P}{\ensuremath{\mathbb{P}}}

\renewcommand{\i}{\mathrm{i}\;}
\def\im{\text{Im\kern1.0pt }}
\def\re{\text{Re\kern1.0pt }}
\def\LMP{Lett. Math. Phys. }
\def\CMP{Commun. Math. Phys. }

\newcommand{\refE}[1]{(\ref{E:#1})}

\newcommand{\refT}[1]{Theorem~\ref{T:#1}}
\newcommand{\refP}[1]{Proposition~\ref{P:#1}}

\newcommand{\btq}{Berezin-Toeplitz quantization}
\newcommand{\w}{\omega}
\newcommand{\Om}{\Omega}
\newcommand{\Ome}{\Omega_\epsilon}
\newcommand{\Ls}{\mathrm{L}^2(M,L)}
\newcommand{\ghm}[1][m]{\Gamma_{hol}(M,L^{#1})}
\newcommand{\gh}{\Gamma_{hol}(M,L)}
\newcommand{\egh}{\mathrm{End}(\Gamma_{hol}(M,L))}
\newcommand{\eghm}{\mathrm{End}(\Gamma_{hol}(M,L^{(m)}))}

\newcommand{\gul}{\Gamma_{\infty}(M,L)}
\newcommand{\gulm}[1][m]{\Gamma_{\infty}(M,L^{#1})}
\newcommand{\cim}{C^{\infty}(M)}
\newcommand{\Tfm}[1][f]{T^{(m)}_{#1}}
\newcommand{\Bfm}[1][f]{B^{(m)}(#1)}
\renewcommand{\d}{\partial}
\newcommand{\db}{\overline{\partial}}

\newcommand{\pnc}[1][N]{\ensuremath{\P^{#1}(\C)}}

\newcommand{\skp}[2]{{\langle #1,#2\rangle}}
\newcommand{\skps}[3]{{\langle #1,#2\rangle}_{#3}}

\begin{document}
\hspace*{\fill} Mannheimer Manuskripte 238

\hspace*{\fill} math.QA/9902066

\vspace*{1.5cm}
\noindent
{ \bf BEREZIN-TOEPLITZ QUANTIZATION AND BEREZIN SYMBOLS}
\noindent
{\bf FOR ARBITRARY COMPACT K\"AHLER MANIFOLDS}
\vspace{1.3cm}\\
\noindent
\hspace*{1in}
\begin{minipage}{13cm}
Martin Schlichenmaier \vspace{0.3cm}\\
Department of Mathematics and Computer Science\\
University of Mannheim, D7, 27\\
D-68131 Mannheim, Germany\\
E-mail:schlichenmaier@math.uni-mannheim.de
\end{minipage}

\vspace*{0.5cm}

\begin{abstract}
\noindent
For phase-space manifolds which are compact K\"ahler 
manifolds relations between the Berezin-Toeplitz quantization 
and the quantization using  Berezin's coherent states and symbols are studied. 
First, results on the Berezin-Toeplitz quantization of
arbitrary quantizable 
compact K\"ahler manifolds due to Bordemann, Meinrenken and
Schlichenmaier are recalled. It is shown that the covariant
symbol map is adjoint to the Toeplitz map.
The Berezin transform for compact K\"ahler manifolds is 
discussed.  

\vspace*{0.2cm}
{\sl  Talk  presented at the XVII${}^{th}$  workshop
on geometric methods in physics,  Bia\l owie\.za, Poland,
July 3 --July 9, 1998
}

\end{abstract}

\section{\hspace{-4mm}.\hspace{2mm}
INTRODUCTION}\label{S:intro}

For phase-space manifolds which are complex K\"ahler manifolds
different quantization schemes of geometric origin
(i.e. related to the complex resp. the K\"ahler structure) have been 
considered. Some of them are connected  with the name of Berezin.
Mainly phase-space manifolds which are either domains in $\C^n$ or 
certain homogeneous spaces were studied. Only recently results on the
quantization of arbitrary  compact K\"ahler manifolds
were obtained. An incomplete list of names related to such results contains
J. Rawnsley, S. Gutt, M. Cahen \cite{Raw},\cite{CGR1},\cite{CGR2},
M. Bordemann, E. Meinrenken, M. Schlichenmaier \cite{BMS},
S. Klimek, A. Lesniewski \cite{KL},
D. Borthwick, A. Uribe \cite{BU}.
For the quantization of additional structures (i.e. the category of
vector bundles over a compact K\"ahler manifold) see also 
E. Hawkins \cite{Haw}.

Here I like to relate for the compact K\"ahler case the Berezin-Toeplitz
quantization (sometimes just called Toeplitz quantization) with the
quantization using Berezin's coherent states, resp. Berezin's
symbols.

Firstly, I will recall results on the \btq\ which were obtained in joint
work with M. Bordemann and E. Meinrenken \cite{BMS}.
I presented them in more detail already 
at the  Bia{\l}owie\.za  1995 workshop \cite{SchlBia95}.
Then I will recall Berezin's coherent states and symbols
\cite{Bere} in their reformulation
and generalization due to Rawnsley \cite{Raw},\cite{CGR1}.
I will show for compact K\"ahler manifolds that
 the Toeplitz operator map and the covariant symbol map
of Berezin are adjoint if one takes the Hilbert-Schmidt norm for the
operators and the, by the epsilon function of 
Rawnsley, deformed Lebesgue measure for the functions.
I will close with introducing  the Berezin transform for 
compact K\"ahler manifolds and discussing certain 
asymptotic properties of it.

In this way I extend results known for the 
bounded symmetric domains 
in $\C^n$ to  arbitrary compact K\"ahler manifolds.  
The study of the Berezin transform 
for such domains in $\C^n$ (hence in some sense at the opposite 
edge of the set of manifolds)
goes back to F. Berezin \cite{Bere} and later  A. Unterberger
and H. Upmeier \cite{UnUp}, M. Engli\v s \cite{Eng},\cite{Eng2} 
 and J. Peetre \cite{EnPe}.

\section{\hspace{-4mm}.\hspace{2mm}
BEREZIN- TOEPLITZ QUANTIZATION}\label{S:btq}

Let $(M,\w)$ be a quantizable K\"ahler manifold. The complex 
manifold $M$ is the phase-space and the K\"ahler form $\w$
(a nondegenerate closed 2-form of type $(1,1)$) is taken as
symplectic form.
Let $(L,h,\nabla)$ be an associated quantum line bundle.
Here $L$ is a holomorphic line bundle, $h$ a hermitian metric
(conjugate linear in the first argument), and $\nabla$ the unique
connection in $L$ which is compatible with the complex structure and the
metric.
With respect to a holomorphic frame of the bundle the metric $h$ can be given 
by a local function $\hat h$ and then the connection will be fixed as
\begin{equation}\label{E:confix}
\nabla_{|}=\d+(\d \log \hat h)+\db\ .
\end{equation}
The quantization condition is the requirement that the curvature form $F$ of 
the line bundle 
coincides with the K\"ahler form of the manifold  up to a factor $(-\i)$, i.e.
\begin{equation}\label{E:quantcond} 
F(X,Y)=\nabla_X\nabla_Y-\nabla_Y\nabla_X-\nabla_{[X,Y]}
=-\i\w(X,Y),
\end{equation}
for arbitrary vector fields $X$ and $Y$.
Again, with respect to a local holomorphic frame this reads as 
\begin{equation}\label{E:wmet}
\w=\i\db\d\log\hat h\ .
\end{equation}
In the following we will restrict our situation  to the compact 
quantizable K\"ahler case.
But note that after some  necessary  modifications
 a lot of the constructions can be extended to
the non-compact case as well.
See \cite{SchlBia95} for more details and for examples.

The first important {\bf observation} is that by the 
quantization condition the line bundle $L$ is a positive line bundle,
resp. in the language of algebraic geometry an ample line bundle.
This says that a positive tensor power $L^{\otimes m_0} $ 
of $L$ is very ample, i.e. $M$ can be embedded as algebraic
submanifold into projective space $\pnc$ 
with the help of a basis $\{s_0,s_1,\ldots,s_N\}$
of the global holomorphic sections of $L^{\otimes m_0} $.
In particular, every quantizable compact K\"ahler manifold is a
projective algebraic manifold. 
Vice versa, every projective manifold is a quantizable K\"ahler manifold with 
K\"ahler form given by  the restriction of 
the K\"ahler form of $\pnc$, (the Fubini-Study form) and  quantum
line bundle given  by the restriction of  the hyperplane bundle.
Again, see \cite{SchlBia95} for details.

Let us assume for the following  that $L$ is  already very ample.
If not, then we  replace $L$ by  $L^{\otimes m_0} $ and the 
K\"ahler form $\w$ by $m_0\,\w$ without changing the complex manifold structure
on $M$. 
Take $\Om=\frac {1}{n!}\w^{\wedge n}$ 
($n=\dim_{\C}M$) as volume form on $M$.
Let $\gul$ be the space of differentiable global sections of $L$ with
scalar product and norm
 \begin{equation}\label{E:skp}
\langle\varphi,\psi\rangle:=\int_M h(\varphi,\psi)\;\Omega\  ,
\qquad
||\varphi||:=\sqrt{\langle \varphi,\varphi\rangle}\ .
\end{equation}
Denote by $\Ls$ the $\mathrm{L}^2$-completion of $\gul$ and by $\gh$ the 
finite-dimensional closed subspace of global holomorphic sections. Let
 \begin{equation}\label{E:proj}
\Pi:\Ls\to\gh,
\end{equation}
be the projection on this subspace.
As usual let $\cim$ be the algebra of complex-valued $C^\infty$-functions on $M$.
Recall that $(\cim,\cdot,\{.,.\})$ is via the symplectic form
$\w$ a Poisson algebra. Its Lie structure $\{.,.\}$
is defined as
\begin{equation}
\{f,g\}:=\w(X_f,X_g),\quad\text{with}\quad \w(X_f,.)=df(.)\ .  
\end{equation}
\begin{definition}\label{D:toeplitz}
(a) For $f\in\cim$ the {\em Toeplitz operator} $T_f$ is defined as
\begin{equation}\label{E:toepop}
T_f:=\Pi\circ(f\cdot):\quad\gh\to\gh,\qquad
s\mapsto \Pi(f\cdot s)\ .
\end{equation}
(b) The map
\begin{equation}\label{E:toepmap}
T:\cim\to\egh,\qquad f\mapsto T_f,
\end{equation}
is called the {\em Berezin-Toeplitz quantization map}.
\end{definition}
In words: The Toeplitz operator $T_f$ associated to the 
differentiable function $f$ 
multiplies the holomorphic section with this function 
and projects the obtained differentiable section back 
to a holomorphic one.

Clearly the map $T$ is linear but it is  neither a 
homomorphism of associative algebras (in general $T_f\cdot T_g\ne T_{fg}$),
nor a Lie algebra homomorphism (in general $T_{\{f,g\}}\ne [T_f,T_g]$).
By the quantization process a lot of  classical informations get lost.
The algebra $\egh$
is finite-dimensional in contrast to the infinite-dimensional algebra  $\cim$.
The right thing to do is to consider instead of only $L$ all its
tensor powers $L^{m}:=L^{\otimes m}$ for $m\in\N$.
With respect to the induced metric $h^{(m)}$ on $L^m$ we obtain now
the corresponding scalar products on the space of global sections of $L^m$,
the projection operator
\begin{equation}
\Pi^{(m)}:\gulm\to\ghm,  
\end{equation}
the Toeplitz operators 
\begin{equation}
\Tfm:\ghm\to\ghm,
\end{equation}
and the Berezin-Toeplitz quantization map
\begin{equation}
T^{(m)}:\cim\to\eghm\ ,
\end{equation}
for every $m\in\N$.
These maps have the correct semi-classical behavior for $m\to\infty$ 
(resp. for $\hbar=\frac 1m\to 0$) as shown by Bordemann, Meinrenken
and Schlichenmaier \cite{BMS}:
\begin{theorem}\label{T:toeapp}
For $f,g\in\cim$ we have 
\begin{alignat}{2}
\mathrm{(a)}\qquad&
\qquad\qquad\qquad&||\Tfm||\qquad\to||f||_\infty,\qquad &m\to\infty,\\
\mathrm{(b)}\qquad&\qquad&||m\;\i[\Tfm,\Tfm[g]]-\Tfm[\{f,g\}]||\to \quad0
\quad,\qquad
&m\to\infty,\\
\mathrm{(c)}\qquad&\qquad&||\Tfm \Tfm[g]-\Tfm[f\cdot g]||\quad\to \quad0\quad
,\qquad &m\to\infty,
\end{alignat}
where for the operators the operator norm and for
the functions the sup-norm have been chosen.
\end{theorem}
For more detailed results on the asymptotics, see 
\cite{BMS},\cite{SchlBia95}.
The proofs employ the theory of generalized 
Toeplitz operators due to 
Boutet de Monvel and Guillemin \cite{BG}.
It is possible to construct a star product (i.e. a 
deformation quantization) with these techniques 
\cite{SchlBia95},\cite{habil}.

Using similar techniques Borthwick and Uribe \cite{BU}
were able to prove the same kind of semi-classical behavior in the
setting of compact symplectic manifolds with almost-complex 
structure.

\refT{toeapp} together with \refP{toesur} implies that this quantization 
is a strict quantization. See for example the recent book by 
N.P. Landsman \cite{Land} for its definition.

\section{\hspace{-4mm}.\hspace{2mm}
BEREZIN COHERENT STATES AND SYMBOLS}\label{S:bcs}

We use the definition of Berezin's coherent states in its 
coordinate independent version and extension 
due to Rawnsley \cite{Raw},\cite{CGR1},
see also \cite{BHSS},\cite{habil}.
Let the situation be as above. In particular, we assume that $L$ is already
very ample. Later we will consider again any $m^{th}$ power of $L$.
Let $\pi:L\to M$ be the bundle projection and $L_0$ the total space of 
$L$ with the zero section $0(M)$ removed.
For $q\in L_0$ fixed and $s\in\gh$ arbitrary we obtain via
 \begin{equation}\label{E:qhdef}
s(\pi(q))=\hat q(s)\cdot q \
\end{equation}
a linear form 
\begin{equation}
\hat q\ :\ \gh\ \to\ \C\ ,\qquad
s\mapsto\hat q(s)\ .
\end{equation}
Hence, with  the scalar product  there exists exactly one holomorphic section 
$e_q$ with 
\begin{equation}\label{E:qhdefa}
\langle e_q,s\rangle=\hat q(s),\qquad\text{for all}\quad s\in\gh\ .
\end{equation}
One calculates
\begin{equation}\label{E:cohtrans}
e_{cq}=\bar c^{-1}\cdot e_q,\qquad c\in\C^*\ .
\end{equation}
The element $e_q$ is called {\em coherent vector (associated to $q$)}.
Note that $e_q\equiv 0$ would imply $\hat q=0$. This yields via \refE{qhdef}
that all sections will vanish at the point $x=\pi(q)$.
But this is a contradiction to the very-ampleness of $L$.
Hence, $e_q\not\equiv 0$ and due to \refE{cohtrans} the element 
$[e_q]:=\{s\in\gh\mid \exists c\in\C^*:s=c\cdot e_q\}$ with $\pi(q)=x$ 
is a well-defined
element of the projective space $\P(\gh)$ only depending on $x\in M$.
It is called the {\em coherent state (associated to $x\in M$)}.

The {\em coherent state embedding} is the 
antiholomorphic embedding
\begin{equation}\label{E:cohemb}
M\quad \to\quad \P(\gh)\ \cong\ \pnc[N],
\qquad 
x\mapsto [e_{\pi^{-1}(x)}].
\end{equation}
In abuse of notation we will understand
in this context under $\pi^{-1}(x)$ always a non-zero element of the 
fiber over $x$.
The coherent state embedding is up to conjugation the Kodaira embedding
with respect to  an orthonormal basis of the sections.
See \cite{BerSchl}  for further considerations of the geometry involved.

We need also the \emph{coherent projectors} used by Rawnsley
\begin{equation}\label{E:cohproj}
P_{\pi(q)}=\frac {|e_q\rangle\langle e_q|}{\langle e_q,e_q\rangle}\ .
\end{equation}
Here we used the convenient bra-ket notation.
{}From \refE{cohtrans} it follows that the projector is well-defined on $M$.
Let us relate the projectors to the metric in the bundle with the help of 
Rawnsley's epsilon function as defined in \cite{Raw}
\begin{equation}\label{E:epsilon}
\epsilon(\pi(q)):=|q|^2\langle e_q,e_q\rangle ,\quad
\text{with}\quad
|q|^2:=h(\pi(q))(q,q).
\end{equation}
Take two sections  $s_1$ and $s_2$. At a fixed point 
$x=\pi(q)$ we can write $s_1(x)=\hat q(s_1)q$ and 
$s_2(x)=\hat q(s_2)q$ and hence for the metric (using \refE{qhdefa})
\begin{equation}\label{E:me}
h(s_1,s_2)(x)=
\overline{\hat q(s_1)}\cdot \hat q(s_2)\cdot |q|^2
=
\langle s_1,e_q\rangle \langle e_q,s_2\rangle |q|^2=
\langle s_1,P_xs_2\rangle \cdot \epsilon(x)\ .
\end{equation}
The {\em covariant Berezin symbol} $\sigma(A)$
of an operator $A\in\egh$ is defined as
\begin{equation}\label{E:covB}
\sigma(A):M\to\C,\qquad
x\mapsto \sigma(A)(x):=Tr(AP_x)=
\frac {\skp {e_q}{Ae_q}}{\skp {e_q}{e_q}},\quad q\in\pi^{-1}(x), q\ne 0\ .
\end{equation}
The symbol $\sigma(A)$ is real-analytic and obeys
\begin{equation}\label{E:cocg}
\sigma(A^*)=\overline{\sigma(A)}\ .
\end{equation}
A closer inspection shows that the linear symbol map
\begin{equation}\label{E:symmap}
\sigma:\egh\to \cim
\end{equation}
is injective, see \cite{CGR1}, or \cite[Prop. 4.1]{BMS}.

Let us now introduce the modified measure 
(note that $\epsilon(x)>0$, for all $x\in M$)
\begin{equation}\label{E:omeps}
\Om_\epsilon:=\epsilon(x)\Om(x)
\end{equation}
on $M$.
The corresponding scalar product on $\cim$ is denoted by ${\skp {.}{.}}_\epsilon$. 
The  {\em contravariant Berezin symbol} $\check \sigma(A)\in\cim$
of an operator is defined by the representation of the operator
$A$ as integral
\begin{equation}\label{E:conB}
A=\int_M\check \sigma(A)(x)P_x\,\Ome(x),
\end{equation}
if such a representation exists.
\begin{proposition}\label{P:contoe}
The Toeplitz operator $T_f$ admits a representation \refE{conB} with
\begin{equation}\label{E:contoe}
\check \sigma(T_f)=f\ , 
\end{equation}
i.e. the function $f$ is the contravariant symbol of the Toeplitz operator
$T_f$.
\end{proposition}
\begin{proof}
Set
\begin{equation}
A:=\int_M f(x)P_x\,\Ome(x)
\end{equation}
then $\check \sigma (A)=f$.
For arbitrary $s_1,s_2\in\gh$ we calculate
(using \refE{me})
\begin{multline}
\skp {s_1}{As_2}=\int_M f(x)\skp {s_1}{P_xs_2}\Ome(x)
=\int_M f(x)h(s_1,s_2)(x)\Om(x)\\
=\int_M h(s_1,fs_2)(x)\Om(x)=
\skp {s_1}{fs_2}=\skp {s_1}{T_fs_2}\ .
\end{multline}
Hence $T_f=A$.
\end{proof}
We introduce on $\egh$ the Hilbert-Schmidt norm 
\begin{equation}\label{E:HS}
\skps {A}{C}{HS}=Tr(A^*\cdot C)\ .
\end{equation}
\begin{theorem}\label{T:adj}
The Toeplitz map $f\to T_f$ and the covariant symbol map
$A\to\sigma(A)$ are adjoint:
\begin{equation}\label{E:adjoint}
\skps {A}{T_f}{HS}=\skps {\sigma(A)} {f} {\epsilon} \ .
\end{equation}
\end{theorem}
\begin{proof}
\begin{equation}
\langle A,T_f\rangle =Tr( A^*\cdot T_f)
=Tr( A^*\int_Mf(x)P_x\,\Omega_\epsilon(x))
=\int_Mf(x) Tr(A^*\cdot P_x)\Omega_\epsilon(x).
\end{equation}
Now applying the definition \refE{covB} and equation \refE{cocg}
\begin{equation}
\langle A,T_f\rangle =\int_M f(x)\sigma( A^*)\Omega_\epsilon(x)
=\int_M \overline{\sigma(A)}(x)f(x)\Omega_\epsilon(x)
={\langle \sigma(A),f(x)\rangle} _\epsilon
\ .
\end{equation}
\end{proof}
The same is valid for every operator $C$ which admits a 
contravariant symbol $\check\sigma(C)$
\begin{equation}
\skps {A}{C}{HS}=\skps {\sigma(A)} {\check\sigma(C)} {\epsilon} \ . 
\end{equation}
In the compact K\"ahler case this will not be an additional result due to the
following
proposition.
\begin{proposition}\label{P:toesur}
The Toeplitz map $f\to T_f$ is surjective. In particular, every operator in 
$\egh$ has a contravariant symbol.
\end{proposition}
\begin{proof}
Choose $A$ an operator orthogonal to $\im T$, i.e. 
$\skp {A}{T_f}=0$ for all $f\in\cim$.
Hence, \refT{adj} implies 
${\skp {\sigma(A)}{f}}_\epsilon=0$  for all $f\in\cim$, i.e. 
$\sigma(A)=0$. By the injectivity of the symbol map this implies $A=0$.
Hence the $T_f$ span the whole $\egh$.
\end{proof}
Now the question  arises when the measure $\Ome$ will be the standard
measure (up to a scalar).
For $M=\pnc$ from the homogeneity of the bundle, and all the
other data it follows 
$\epsilon\equiv const$. By a result of Rawnsley \cite{Raw},
 resp. Cahen, Gutt and
Rawnsley \cite{CGR1},
$\epsilon\equiv const$ if and only if the quantization  is 
projectively induced. 
This means that using the conjugate of the coherent state embedding, 
the K\"ahler form $\w$ of $M$ coincides with the pull-back of the
Fubini-Study form. Note that not every quantization is
projectively induced, see the discussion in \cite{SchlBia95}.

{\bf Appendix.}
To compare the global description with  Berezin's original approach
we have to choose a section $s_0\in\gh$, $s_0\not\equiv 0$.
On the open set $\ V:=\{x\in M\mid s_0(x)\ne 0\}$.
the section $s_0$ is a holomorphic frame for the bundle $L$.
Hence, every holomorphic (resp. differentiable) section $s$ can be described  
as $s(x)=\hat s(x)s_0(x)$ with a holomorphic (resp. differentiable) 
function on $V$. The map $s\mapsto \hat s$ gives an isometry
of $\gh$, resp. $\gul$ with the L${}^2$ space of holomorphic, resp.
differentiable functions on $V$ with respect to the measure
\begin{equation}
 \mu_{s_0}(x)=h(s_0,s_0)\Omega(x)\ . 
\end{equation}
The scalar product can be given as
\begin{equation}
\skp {s}{t}=\int_V\overline{\hat s}\cdot{\hat t}\cdot h(s_0,s_0)\Omega(x)=
\int_V\overline{\hat s}\cdot {\hat t}\cdot\exp(-K(x))\Omega(x)\ .
\end{equation}
Here $K$ is a local K\"ahler potential with respect to the frame $s_0$.
It is given by 
\begin{equation}
K(x)=-\log h(s_0,s_0)(x)\ . 
\end{equation}
Note that by the quantization condition \refE{wmet} we have indeed 
$\w=\i\d\db K$.
\section{\hspace{-4mm}.\hspace{2mm}
THE BEREZIN TRANSFORM}\label{S:btrans}

Starting from $f\in\cim$ we can assign to $f$ its Toeplitz operator
$T_f\in\egh$ and then assign to $T_f$ the covariant symbol $\sigma(T_f)$ which 
is again an element of $\cim$.  
Altogether we obtain a map $f\mapsto B(f)=\sigma(T_f)$.
This map is called {\em Berezin transform}.

Recall that $f$ is the contravariant symbol of the Toeplitz operator $T_f$.
Hence $B$ gives a correspondence between contravariant symbols
and covariant symbols of operators.
The Berezin transform was introduced and studied by 
Berezin \cite{Bere} for certain classical symmetric 
domains in $\C^n$. These results where 
extended by Unterberger and Upmeier \cite{UnUp},
see also Engli\v s \cite{Eng},\cite{Eng2} and Engli\v s and Peetre \cite{EnPe}.
As seen above the Berezin transform makes sense also
in the compact K\"ahler case which we
consider here.
Let me study it in some more detail.

For $s,t$ holomorphic sections of $L$  and $f\in\cim$
we have 
$\skp {t}{T_fs}=\skp {t}{fs}$. Hence,
\begin{equation}\label{E:btr1}
B(f)(x)=\sigma(T_f)(x)=
\frac {\skp {e_q}{fe_q}}{\skp {e_q}{e_q}},\quad x=\pi(q)\ , q\ne 0\ .
\end{equation}
It is possible to find another useful description.
{}From \refP{contoe}  follows
\begin{equation}
T_f=
\int_M f(y) \frac {|e_{q'}\rangle\langle e_{q'}|}{\langle e_{q'},e_{q'}\rangle}
\,\Ome(y) ,\quad y=\pi(q')\ , q'\ne 0\ ,
\end{equation}
and we obtain with $q\in\pi^{-1}(x), q\ne 0$ 
\begin{equation}\label{E:btr2}
B(f)(x)= \frac {\skp {e_{q}}{T_fe_{q}}}{\skp {e_{q}}{e_{q}}}
=\int_M f(y) \frac {\skp {e_{q'}}{e_{q}}}{\skp {e_{q'}}{e_{q'}}}
\frac {\skp {e_{q}}{e_{q'}}}{\skp {e_{q}}{e_{q}}}
\,\Ome(y) \ .
\end{equation}
We set 
\begin{equation}
\label{E:btkern}
K(x,y):=\frac  {{|\skp {e_{q'}}{e_q}|}^2}{\skp {e_{q'}}{e_{q'}}
\skp {e_{q}}{e_{q}}}\ .
\end{equation}
Clearly, $0\le K(x,y)\le 1$ for every $x,y\in M$.
Note that $K$ is one of  the two-point functions 
introduced in \cite{CGR2}.
Now the Berezin transform can be written 
with the help of an integral kernel  as
\begin{equation}
\label{E:btint}
B(f)(x)=\int_M f(y) K(x,y)\Ome(y)\ .
\end{equation}
Everything can be done for any positive power $m$ of the
line bundle $L$.
\begin{definition}
Let $m\in\N$. The map
\begin{equation}\label{E:btrans}
B:\cim\to\cim,\qquad f\mapsto B^{(m)}(f)=\sigma^{(m)}(\Tfm)
\end{equation}
is called the {\em Berezin transform of level} $m$.
\end{definition}
Note that (nearly) everything depends on $m$,
the scalar product, the coherent states $e_q^{(m)}$, the 
symbol maps, the epsilon function,  the integral
kernel $K^{(m)}(x,y)$, etc.
The asymptotic behavior of 
$B^{(m)}$, resp. of $B^{(m)}(f)$ as $m\to\infty$
contains interesting information.
In the classical bounded symmetric  domain case in $\C^n$
considered
by Berezin \cite{Bere} the asymptotics contains a lot of information
about the quantization.
There the first two terms in the asymptotics are corresponding to the
fact, that the quantization has the correct semi-classical behavior.
See also Engli\v s \cite{Eng2} for results on  some more general
domains in $\C^n$.

In the compact K\"ahler case we have
\begin{theorem}\label{T:symbols}
Let $f\in\cim$ then 
\begin{equation}\label{E:symine}
|\Bfm|_\infty=|\sigma^{(m)}(\Tfm)|_\infty\quad\le \quad||\Tfm||\quad\le\quad   
|\check\sigma(\Tfm)|_\infty=|f|_\infty\ .
\end{equation}
\end{theorem}
\begin{proof}
The two out-most equalities are by definition.
To simplify the notation let us drop the super-script $(m)$ in the proof.
Using the Cauchy-Schwarz inequality we calculate ($x=\pi(q)$)
\begin{equation}
| \sigma(T_f)(x)|^2=
\frac {|\skp {e_q}{T_fe_q}|^2}{{\skp {e_q}{e_q}}^2}\le
\frac {\skp {T_fe_q}{T_fe_q}}{\skp {e_q}{e_q}}\le
||T_f||^2\ .
\end{equation}
Here the last step follows from the definition of the operator norm
\begin{equation}
||T_f||=\sup_{\scriptstyle s\in\gh\atop \scriptstyle s\ne 0}
\frac {||T_fs||}{||s||}\ .  
\end{equation}
This shows the first inequality in \refE{symine}.
For the second inequality introduce the multiplication
operator $M_f$ on $\gul$. Then 
$\ ||T_f||=||\Pi\,M_f\,\Pi||\le ||M_f||\ $ and
for  $\varphi\in\gul$,  $\varphi\ne 0$
\begin{equation}
\frac {{||M_f \varphi||}^2}{||\varphi||^2}=
\frac {\int_M h(f \varphi,f \varphi)\Omega}
 {\int_M h(\varphi,\varphi)\Omega}
=
\frac {\int_M \overline{f(z)}f(z)h(\varphi,\varphi)\Omega}
{\int_M h(\varphi,\varphi)\Omega}
\le
|f|{}_\infty^2\ .  
\end{equation}
Hence,
\begin{equation}
||T_f||\le ||M_f||=\sup_{\varphi\ne 0}
\frac {||M_f\varphi||}{||\varphi||}\le |f|_\infty .
\end{equation}
\end{proof}
\begin{corollary}
Let $A$ be an operator of $\ghm$ then
\begin{equation}
|\sigma^{(m)}(A)|_\infty\quad\le\quad ||A||\quad\le \quad
 |\check\sigma^{(m)}(A)|_\infty\ . 
\end{equation}
\begin{proof}
By \refP{toesur}   
every operator is a Toeplitz operator. Hence we can apply \refT{symbols}.
\end{proof}
\end{corollary}
{\em  Warning:}
The statement of \refP{toesur} should not be misinterpreted 
in the way that given
a (natural) family of operators $A^{(m)}$ we will have $A^{(m)}=\Tfm$.
It only states that for a fixed $m$ there is a function $f^{(m)}$ such that 
$A^{(m)}=T_{f^{(m)}}^{(m)}$.
An important example is given by the operator $Q^{(m)}$ of geometric quantization
(with K\"ahler polarisation). By a result of Tuynman we have 
$\ Q^{(m)}=T^{(m)}_{\i(f-\frac 1{2m}\Delta f)}$ \cite{Tuyn}, see also
\cite{BHSS}.

By the above results  we see that 
$B^{(m)}$ is norm contracting.
 In particular $\Bfm$ is bounded  by $|f|_\infty$.
If we consider the asymptotic expansion 
\begin{equation}
  \label{E:asymp}
\Bfm=A_0(f)+A_1(f)\frac 1m+A_2(f)\frac {1}{m^2}+\cdots,  
\end{equation}
then in the 
non-compact K\"ahler case of bounded symmetric domains \cite{UnUp}
and for the planar domains of hyperbolic type in $\C$ \cite{Eng}
 it is shown  that $A_0=id$ and that the $A_i$ are
polynomials in the 
 invariant differential operators (generalized
Laplacians), resp. in the Laplace-Beltrami operator 
only depending on the geometry of $M$.
In particular note, that the fundamental regions 
for compact Riemann surfaces of genus $g\ge 2$ are 
planar domains of hyperbolic type. 
Similar results are expected also for the general compact K\"ahler case.
Let me quote one result. For $\epsilon\equiv const$, i.e. the
projectively 
induced quantization case, all $\epsilon^{(m)}$ will (individually) be constant.
{}From \refE{btr1} it follows with the help of the stationary
phase theorem in a way similar to the proof of part (a) of \refT{toeapp}
(see \cite{BMS}, \cite[p.83]{habil}, \cite{SchlBia95})
\begin{equation}
  \label{E:Bas}
  (\Bfm)(x)=f(x)+O(\frac 1m)\ .
\end{equation}
Clearly, from an  asymptotics given by 
\refE{Bas} part (a) of \refT{toeapp}  follows with the help of \refT{symbols}.
Hence, as expected the zero order part of the asymptotics of $B^{(m)}$  is related to  
the first condition of the correct semi-classical behavior.
Details and further results have to be postponed to a forthcoming paper.
Let me close with the remark that 
the asymptotic expansion  \refE{asymp} with $A_0=id$  
gives a relation between  the star product obtained via
Berezin-Toeplitz quantization and the star product obtained via
Berezins's recipe \cite{Bere},\cite{CGR2}. 
See also \cite{CGR2} for results partly related to the asymptotics.


\end{document}